\documentclass[12pt, a4paper]{amsart}
\usepackage{amsmath}
\usepackage{amssymb}
\usepackage{amsthm}
\usepackage{amscd}
\usepackage[all]{xy}

\theoremstyle{plain}
\newtheorem{thm}{Theorem}[section]
\newtheorem{prop}[thm]{Proposition}

\newtheorem{ques}[thm]{Question}
\newtheorem{mainthm}{Theorem}

\theoremstyle{definition}
\newtheorem{defn}[thm]{Definition}
\newtheorem{expl}[thm]{Example}
\theoremstyle{remark}
\newtheorem{rem}[thm]{Remark}

\newtheorem*{notation}{Notation}

\newtheorem{cln}{Claim}

\newcommand{\mC}{{\mathbb C}}

\newcommand{\mQ}{{\mathbb Q}}
\newcommand{\mR}{{\mathbb R}}

\newcommand{\mZ}{{\mathbb Z}}

\newcommand{\Spec}{\mathrm{Spec}\,}

\newcommand{\Supp}{\mathrm{Supp}\,}

\renewcommand{\labelenumi}{(\arabic{enumi})}
\numberwithin{equation}{section}

\begin{document}
\title{When does the subadditivity theorem\\ for multiplier ideals hold?}
\author{Shunsuke Takagi}
\address{Graduate School of Mathematical Sciences, University of Tokyo, 3-8-1, Komaba, Meguro, Tokyo 153-8914, Japan}
\email{stakagi@ms.u-tokyo.ac.jp}
\author{Kei-ichi Watanabe}
\address{Department of Mathematics, College of Humanities and Sciences, 
Nihon University, Setagaya-ku, Tokyo 156--0045, Japan}
\email{watanabe@math.chs.nihon-u.ac.jp}
\subjclass[2000]{Primary 13B22; Secondary 14J17}
\thanks{Both authors thank MSRI for the support and hospitality during their stay in the fall of 2002. The second author was partially supported by Grants-in-Aid in Scientific Researches, 13440015, 13874006; and his stay at MSRI was supported by the Bunri Fund, Nihon University.}
\baselineskip 15pt
\footskip = 32pt

\begin{abstract}
Demailly, Ein and Lazarsfeld \cite{DEL} proved the subadditivity theorem for multiplier ideals on non-singular varieties, which states the multiplier ideal of the product of ideals is contained in the product of the individual multiplier ideals. We prove that, in two-dimensional case, the subadditivity theorem holds on log terminal singularities. However, in higher dimensional case, we have several counterexamples. We consider the subadditivity theorem for monomial ideals on toric rings and construct a counterexample on a three-dimensional toric ring. 
\end{abstract}

\maketitle
\markboth{SHUNSUKE TAKAGI AND KEI-ICHI WATANABE}{THE SUBADDITIVITY THEOREM FOR MULTIPLIER IDEALS}

\section*{Introduction}
Multiplier ideals were first introduced in the complex analytic context in the work of Demailly, Nadel, Siu and others, and they proved a Kodaira-type vanishing theorem involving these ideals. Multiplier ideals can be reformulated in a purely algebro-geometric setting in terms of resolution of singularities and discrepancy divisors, and nowadays this notion has become a fundamental tool in birational geometry.

Demailly, Ein and Lazarsfeld \cite{DEL} proved the subadditivity theorem for multiplier ideals on non-singular varieties, which states the multiplier ideal of the product of ideals is contained in the product of the individual multiplier ideals. This theorem itself is miraculous for commutative algebraists, and moreover it has several interesting applications to commutative algebra and algebraic geometry. For example, the problem concerning the growth of symbolic powers of ideals in regular local rings (see \cite{ELS}), Fujita's approximation theorem which asserts that most of the volume of a big divisor can be accounted for by the volume of an ample $\mQ$-divisor on a modification (see \cite{Fu} and \cite{La}), etc.
However, their proof of the subadditivity theorem works only on non-singular varieties over a field of characteristic zero, because their proof needs the Kawamata-Viehweg vanishing theorem and the fact that the diagonal embedding is a complete intersection. 
Hence we investigate when the subadditivity theorem holds on singular varieties which admit a resolution of singularities. 
The multiplier ideal associated to the unit ideal defines the locus of non-log-terminal points. 
Therefore, on non-log-terminal singularities, the subadditivity theorem fails.
Conversely, in two-dimensional case, using a characterization of integrally closed ideals via anti-nef cycles, we show that the subadditivity theorem holds on log terminal singularities which are not necessarily essentially of finite type over a field of characteristic zero. 

\renewcommand{\themainthm}{2.2}
\begin{mainthm}
Let $(A, \mathfrak{m})$ be a two-dimensional $\mQ$-Gorenstein normal local ring. Then $A$ is log terminal if and only if the subadditivity theorem holds, that is, for any two ideals $\mathfrak{a}$, $\mathfrak{b} \subseteq A$, 
$$\mathcal{J}(\mathfrak{ab}) \subseteq \mathcal{J}(\mathfrak{a}) \mathcal{J}(\mathfrak{b}).$$
\end{mainthm}

However, in higher dimensional case, we have several counterexamples to Theorem \ref{m-primary} (See Example \ref{3-dim}). So we investigate the subadditivity theorem for monomial ideals. 
The multiplier ideal associated to a monomial ideal is characterized by the Newton polygon (see \cite{HY} and \cite{How}) and it is easy to calculate this ideal. We expected that the subadditivity theorem for monomial ideals might hold on all toric rings. But, unfortunately, we found a counterexample on a three-dimensional toric ring (see Example \ref{toric}). 

\section{Multiplier ideals}
\begin{notation}
Throughout this paper, let $(A, \mathfrak{m})$ be an excellent normal $\mQ$-Gorenstein local ring satisfying one of the following conditions:
\begin{itemize}
\item $(A, \mathfrak{m})$ is two-dimensional.
\item $(A, \mathfrak{m})$ is essentially of finite type over a field of characteristic zero.
\end{itemize}
\end{notation}
First we recall the definition of multiplier ideals. Refer to \cite{La} for the general theory of multiplier ideals.
\begin{defn}
Let $\mathfrak{a}$ be an ideal in $A$. By \cite{Hi}, \cite{Li1} and \cite{Li2}, there exists a resolution of singularities $f:X \to \Spec A$ such that the ideal sheaf $\mathfrak{a} \mathcal{O}_X=\mathcal{O}_X(-F)$ is invertible and $\mathrm{Exc}(f) \cup \Supp F$ is a simple normal crossing divisor, where $\mathrm{Exc}(f)$ is the exceptional locus of $f$.
Fix a rational number $c > 0$.
Then the \textit{multiplier ideal}\footnote{Lipman \cite{Li3} calls this ideal the ``adjoint ideal.'' However, Lazarsfeld \cite{La} uses the term ``adjoint ideal'' in a different sense.  To avoid confusion, we adopt the term ``multiplier ideal'' in this paper.} associated to $c$ and $\mathfrak{a}$ is defined to be
$$\mathcal{J}(\mathfrak{a}^c) = \mathcal{J}(A, \mathfrak{a}^c)=H^{0}(X,\mathcal{O}_X(\lceil K_X -f^*K_A -cF \rceil)) \subseteq A,$$
where $K_X$ and $K_A$ are the canonical divisors of $X$ and $\Spec A$ respectively.
In particular, $A$ is said to be a {\it log terminal singularity} if $\mathcal{J}(A) =A$.

Similarly we can also define the multiplier ideal $\mathcal{J}( \mathfrak{a}^c \mathfrak{b}^d )$ associated to two ideals $\mathfrak{a}$, $\mathfrak{b}$ in $A$ and two rational numbers $c$, $d>0$: let $f:X \to \Spec A$ be a resolution of singularities such that $\mathfrak{a}\mathcal{O}_X=\mathcal{O}_X(-F_a)$ and $\mathfrak{b}\mathcal{O}_X=\mathcal{O}_X(-F_b)$ are invertible and $\mathrm{Exc}(f) \cup \Supp F_a \cup \Supp F_b$ is a simple normal crossing divisor. Then 
$$\mathcal{J}( \mathfrak{a}^c \mathfrak{b}^d )=H^{0}(X,\mathcal{O}_X(\lceil K_X -f^*K_A -cF_a-dF_b \rceil)) \subseteq A.$$

\end{defn}

\begin{rem}
\begin{enumerate}
\item Multiplier ideals are independent of the choice of a desingularization $f:X \to \Spec A$.
\item Log terminal singularities are rational singularities.
\end{enumerate}
\end{rem}
The following basic properties of multiplier ideals immediately follow.
\begin{prop}\label{basic}
Let $\mathfrak{a}$ and $\mathfrak{b}$ be ideals in $A$, and $c>0$ be a rational number.
\begin{enumerate}
\renewcommand{\labelenumi}{(\roman{enumi})}
\item If $\mathfrak{a} \subseteq \mathfrak{b}$, then $\mathcal{J}(\mathfrak{a}^c) \subseteq \mathcal{J}(\mathfrak{b}^c)$. 

\item $\mathcal{J}(\mathfrak{a}^c)$ is integrally closed. Moreover $\mathcal{J}(\mathfrak{a}^c)=\mathcal{J}(\overline{\mathfrak{a}}^c)$, where we denote by $\overline{\mathfrak{a}}$ the integral closure of $\mathfrak{a}$.

\item Suppose that $A$ is a log terminal singularity. Then $\mathfrak{a} \subseteq \mathcal{J}(\mathfrak{a})$. Furthermore, if $\mathfrak{a}$ is an ideal of pure height one, then $\mathcal{J}(\mathfrak{a})=\mathfrak{a}$.

\end{enumerate}
\end{prop}
\begin{proof}
We will show only (iii). Let $f:X \to \Spec A$ be a resolution of singularities such that $\mathfrak{a}\mathcal{O}_X=\mathcal{O}_X(-Z)$ is invertible and $\mathrm{Exc}(f) \cup \Supp Z$ is a simple normal crossing divisor.
Since $A$ is log terminal, 
$$\mathcal{J}(\mathfrak{a})= H^0(X, \mathcal{O}_X(\lceil K_X-f^*K_A - Z \rceil))
                         \supseteq H^0(X, \mathcal{O}_X(-Z))
                         =\bar{\mathfrak{a}}.$$
Since $\mathrm{codim}_A (\Supp (\mathcal{J}(\mathfrak{a})/\bar{\mathfrak{a}})) \ge 2$, we have $\mathfrak{a}=\bar{\mathfrak{a}}=\mathcal{J}(\mathfrak{a})$ when $\mathfrak{a}$ is divisorial.
\end{proof}

Demailly, Ein and Lazarsfeld proved the following theorem, which is called the subadditivity theorem.
\begin{thm}[\cite{DEL}]\label{DEL}
Let $A$ be a regular local ring essentially of finite type over a field of characteristic zero, and let $\mathfrak{a}$ and $\mathfrak{b}$ be any two ideals in $A$. Fix any rational numbers $c$, $d>0$. Then
$$\mathcal{J}(\mathfrak{a}^c \mathfrak{b}^d) \subseteq \mathcal{J}(\mathfrak{a}^c)  \mathcal{J}(\mathfrak{b}^d).$$
\end{thm}
\begin{rem}
Demailly, Ein and Lazarsfeld use the Kawamata-Viehweg vanishing theorem, hence the condition that $A$ is essentially of finite type over a field of characteristic zero is necessary for their proof.
\end{rem}
In this paper, we say that {\it the  subadditivity theorem holds} if $\mathcal{J}(\mathfrak{a} \mathfrak{b}) \subseteq \mathcal{J}(\mathfrak{a})  \mathcal{J}(\mathfrak{b})$ for any ideals $\mathfrak{a}, \mathfrak{b} \subseteq A$, and {\it the strong subadditivity theorem holds} if $\mathcal{J}(\mathfrak{a}^c \mathfrak{b}^d) \subseteq \mathcal{J}(\mathfrak{a}^c)  \mathcal{J}(\mathfrak{b}^d)$ for any $\mathfrak{a}, \mathfrak{b}$ and any rational numbers $c, d>0$.

\section{two-dimensional case}
In this section, we investigate when the subadditivity theorem holds in two-dimensional case. The following characterization of integrally closed ideals is quite useful.

\begin{thm}[\cite{Li1}, \cite{Gi}]\label{anti-nef}
Let $(A, \mathfrak{m})$ be a two-dimensional rational singularity and fix a resolution of singularities $f:X \to \Spec A$ with $E$ the exceptional locus on X. Let $E_1, \dots, E_r$ be all the irreducible components of $E$. Then there is a one-to-one correspondence between the set of integrally closed ideals $I$ in $A$ such that $I\mathcal{O}_X$ is invertible and the set of effective $f$-anti-nef cycles $Z$ on $X$ $($i.e. $Z \ge 0$ and $Z \cdot E_i \leq 0$ for all $1 \le i \le r$$)$. The correspondence is given by $I\mathcal{O}_X=\mathcal{O}_X(-Z)$ and $I=H^0(X, \mathcal{O}_X(-Z))$.
\end{thm}

Using the above theorem, we obtain necessary and sufficient conditions for the subadditivity theorem to hold.

\begin{thm}\label{m-primary}
Let $(A, \mathfrak{m})$ be a two-dimensional $\mQ$-Gorenstein normal local ring. Then $A$ is a log terminal singularity if and only if the subadditivity theorem holds, that is, for any two ideals $\mathfrak{a}, \mathfrak{b} \subseteq A$,
$$\mathcal{J}(\mathfrak{ab}) \subseteq \mathcal{J}(\mathfrak{a}) \mathcal{J}(\mathfrak{b}).$$
\end{thm}
\begin{proof}
If the subadditivity theorem holds, then $\mathcal{J}(A) \subseteq \mathcal{J}(A)^2$. Thus $\mathcal{J}(A)=A$, namely $A$ is log terminal. 
Hence we will show the converse implication, that is, we will prove that for any two ideals $\mathfrak{a}, \mathfrak{b} \subseteq A$, $\mathcal{J}(\mathfrak{ab})\subseteq \mathcal{J}(\mathfrak{a}) \mathcal{J}(\mathfrak{b})$ when $A$ is a log terminal singularity. 
By Proposition \ref{basic} (ii), we may assume that $\mathfrak{a}$ and $\mathfrak{b}$ are integrally closed.
Let $f:X \to \Spec A$ be a resolution of singularities such that $\mathfrak{a}\mathcal{O}_X=\mathcal{O}_X(-F_a)$ and $\mathfrak{b}\mathcal{O}_X=\mathcal{O}_X(-F_b)$ are invertible and $\mathrm{Exc}(f) \cup \Supp F_a \cup \Supp F_b$ is a simple normal crossing divisor. By Theorem \ref{anti-nef}, $F_a$ and $F_b$ are $f$-anti-nef cycles on $X$, which are not necessarily supported on the exceptional locus of $f$.
By the definition of multiplier ideals, denoting by $K$ the relative canonical divisor $K_{X}-{f}^{*}K_{A}$ of $f$, we have
\begin{align*}
\mathcal{J}(\mathfrak{a}) \mathcal{J}(\mathfrak{b}) &= H^0(X,\mathcal{O}_X(\lceil K \rceil -F_a)) \cdot H^0(X,\mathcal{O}_X(\lceil K \rceil -F_b)), \\
\mathcal{J}(\mathfrak{ab}) &= H^0(X,\mathcal{O}_X(\lceil K \rceil -F_a - F_b)).
\end{align*}
Here, for every cycle $Z$ on $X$, we denote by $\mathrm{an}_f(Z)$ the $f$-anti-nef closure of $Z$, namely the minimal $f$-anti-nef cycle among all cycles on $X$ which is bigger than or equal to $Z$. 
Note that $\mathrm{an}_f(Z)$ is uniquely determined by $Z$ (cf. \cite{Ar}).
Since $A$ is a rational singularity, the product of integrally closed ideals of $A$ is also integrally closed \cite{Li1}. 
Hence $\mathcal{J}(\mathfrak{a}) \mathcal{J}(\mathfrak{b})$ and $\mathcal{J}(\mathfrak{ab})$ are integrally closed, and by Theorem \ref{anti-nef} again, 
$\mathcal{J}(\mathfrak{a}) \mathcal{J}(\mathfrak{b})$ and $\mathcal{J}(\mathfrak{ab})$ correspond to the cycles
$\mathrm{an}_f(F_a-\lceil K \rceil)+\mathrm{an}_f(F_b - \lceil K \rceil)$ and $\mathrm{an}_f(F_a + F_b -\lceil K \rceil)$ respectively. 
Therefore it suffices to show that
\begin{align}
\mathrm{an}_f(F_a-\lceil K \rceil)+\mathrm{an}_f(F_b - \lceil K \rceil) \leq  \mathrm{an}_f(F_a + F_b -\lceil K \rceil).
\end{align}
In order to prove this, we prepare some notations.
We can assume that the residue field $A/\mathfrak{m}$ is algebraically closed.
Then the morphism $f$ can be factorized as follows.
$$X:=X_n \xrightarrow{f_n} X_{n-1} \xrightarrow{f_{n-1}} \cdots \xrightarrow{f_1} X_0 \xrightarrow{f_0} \Spec A,$$
where $f_i:X_i \to X_{i-1}$ is a contraction of a $(-1)$-curve $E_i$ on $X_i$ for every $i=1, \dots, n$ and $f_0 : X_0 \to \Spec A$ is the minimal resolution of $\Spec A$. 
We denote by $\pi_i:X \to X_i$ the composite of $f_{i+1}, \dots, f_n$ for $i=0,1, \dots, n-1$ and by $\pi_{i,j}:X_i \to X_j$ the composite of $f_{j+1}, \dots, f_i$ for $i >j$.
Then the relation $\vartriangleright$ on $f$-exceptional divisors $E_1, \cdots, E_n$ is defined as follows:
$E_i \vartriangleright E_j$ if and only if the intersection number ${\pi_i}^*E_i \cdot {\pi_j}^{-1}_*E_j$ is positive, where ${\pi_j}^{-1}_*E_j$ is the strict transform of $E_j$ on $X$.
Since the relation $\vartriangleright$ is not a order relation, we denote by $>$ the order relation generated by $\vartriangleright$.

Let $\mathbf{P}$ be the proximity matrix, that is, the matrix $\mathbf{P}:=(p_{ij})_{1\le i,j\le n}$ given by 
$$p_{ij}=
\left\{
\begin{array}{cl}
1 & \text{if $i=j$} \\
-1 & \text{if $E_j \vartriangleright E_i$} \\
0 & \text{otherwise}
\end{array}
\right.
$$
(see \cite{DV} and \cite{Li4}).
\begin{cln}
Let $Z$ be a cycle on $X$ and we write
$$Z={\pi_0}^*{\pi_0}_*Z+\sum_{i=1}^n d_i{\pi_i}^*E_i.$$
Then $Z$ is $f$-anti-nef if and only if ${\pi_{0}}_*Z$ is an $f_0$-anti-nef cycle on $X_0$ and every component of $\mathbf{P d}$ is a nonnegative integer, where $\mathbf{d}={}^t(d_1,\dots,d_n)$.
Notice that if $Z$ is $f$-anti-nef, then $d_i \ge 0$ for every $i$.
\end{cln}
\begin{proof}[Proof of Claim $1$]
Since $Z={\pi_{j-1}}^*{\pi_{j-1}}_*Z+\sum_{i=j}^n d_i{\pi_i}^*E_i$ for every $j=1,2, \dots, n$,
$$Z \cdot {\pi_j}^{-1}_*E_j=\sum_{i=j}^n d_i E_i \cdot {\pi_{i,j}}^{-1}_*{E_j}=-\mathbf{e}_j\mathbf{Pd},$$
where $\mathbf{e}_j$ is the $n$-tuple (row) vector such that the $j$-th component is one and other components are zero. 
On the other hand, for each $f_0$-exceptional curve $F$, we have $Z \cdot {\pi_0}^{-1}_*F={\pi_0}_*Z \cdot F$. Therefore $Z$ is $f$-anti-nef if and only if $\mathbf{e}_j\mathbf{Pd}$ is nonnegative for all $j=1,\dots,n$ and ${\pi_0}_*Z$ is $f_0$-anti-nef.
\end{proof}

Fix any effective $f$-anti-nef cycle $Z$ on $X$ and write $Z={\pi_{0}}^*{\pi_0}_*Z+\sum_{i=1}^n d_i{\pi_i}^*E_i$. Then we investigate a process for computation of the anti-nef closure of $Z-\lceil K \rceil$. Let $\Lambda$ be a subset of $\{1,2, \dots, n\}$ such that $\lceil K \rceil=\sum_{i \in \Lambda} {\pi_i}^*E_i$.
Then define $Z^{(k)}:={\pi_0}^*{\pi_0}_*Z+\sum_{i=1}^n d_i^{(k)}{\pi_i}^*E_i$ inductively as follows:
Let 
\begin{align}
d_i^{(0)}=\left\{
\begin{array}{ll}
d_i-1 & \text{if $i \in \Lambda$ and $d_i > 0$} \\
d_i & \text{otherwise}
\end{array}
\right..
\end{align}
for all $i=1, \dots, n$.
If $Z^{(1)}, \cdots, Z^{(k-1)}$ are defined and if $\mathbf{e}_j\mathbf{P}\mathbf{d}^{(k-1)}$ is negative for some $1 \le j \le n$ where $\mathbf{d}^{(k-1)}={}^t(d_1^{(k-1)},\dots,d_n^{(k-1)})$, then choose one of such $j$ (we denote this by the same letter $j$) and set 
\begin{align}
d_i^{(k)}=\left\{
\begin{array}{ll}
d_j^{(k-1)}+1 & \text{if $i=j$} \\
d_i^{(k-1)}-1 & \text{if $E_i \vartriangleright E_j$ and $d_i^{(k-1)} > 0$} \\
d_i^{(k-1)}   & \text{otherwise}
\end{array}
\right..
\end{align}

\begin{cln}
This process stops after finitely many steps. When the process stops at $Z^{(k_0)}$, the cycle $Z^{(k_0)}$ is the $f$-anti-nef closure of $Z-\lceil K \rceil$.
\end{cln}
\begin{proof}[Proof of Claim $2$]
This is similar to the computation of the fundamental cycle (see \cite{Ar}). So we give only one remark here. 
Some readers may think that for the minimality of the anti-nef closure of $Z-\lceil K \rceil$, $d_i^{(0)}$ should be defined as follows: 
\begin{small}
$$d_i^{(0)}=
\left\{
\begin{array}{ll}
d_i-1 & \text{if $i \in \Lambda$ and $d_i>0$}\\
d_i-1 & \text{if $d_i>0$ and if $d_j=0$ for some $j \in \Lambda$ such that $E_i \vartriangleright E_j$} \\
d_i   & \text{otherwise}
\end{array}
\right.$$
\end{small}
However, this definition coincides with the above one.
In fact, since $Z$ is an $f$-anti-nef cycle, once $d_j=0$ we have $d_i=0$ for all $i$ such that $E_i \vartriangleright E_j$. Therefore when $d_i >0$, there exists no such $E_i \vartriangleright E_j$ as $d_j=0$.
\end{proof}

In this paper, we call such a sequence of cycles as $\{Z^{(0)}, \cdots, Z^{(k_0)}\}$ a computation sequence for $Z-\lceil K \rceil$ (there may be several computation sequences).
\begin{cln}
$d_i-1 \le d_i^{(k)} \le d_i$ for every $i=1,\dots,n$ and $k=0,1,\dots,k_0$. 
Therefore $\mathbf{e}_j\mathbf{P}\mathbf{d}^{(k)} <0$ if and only if $\mathbf{e}_j\mathbf{P}\mathbf{d}=0$, $d_j^{(k)}=d_j-1$ and $d_i^{(k)}=d_i$ for all $i$ such that $E_i \vartriangleright E_j$.
\end{cln}
\begin{proof}[Proof of Claim $3$]
Since $Z$ is $f$-anti-nef, $\mathrm{an}_f(Z-\lceil K \rceil) \le Z$ and $d^{(i)} \le d_i$.
We will show that $d_i-1 \le d_i^{(k)}$ by induction on $k$. When $k=0$, the assertion is trivial. Hence we may assume that $k \ge 1$. By the induction hypothesis, $d_i-1 \le d_i^{(k-1)} \le d_i$. 
If $d_i^{(k-1)}=d_i-1$, then for every $E_j \vartriangleleft E_i$, $\mathbf{e}_j\mathbf{P}\mathbf{d}^{(k-1)} \ge 0$. Therefore $d_i^{(k)} \ge d_i-1$.
\end{proof}
Now using this process, we will prove the equation $(2.1)$. 
Write 

\begin{align*}
F_a&={\pi_0}^*{\pi_0}_*F_a+\sum_{i=1}^{n}a_i{\pi_i}^*E_i, \\
F_b&={\pi_0}^*{\pi_0}_*F_b+\sum_{i=1}^{n}b_i{\pi_i}^*E_i,
\end{align*}
and denote $\mathbf{a}={}^t(a_1,\dots,a_n)$ and $\mathbf{b}={}^t(b_1,\dots,b_n)$.

Let $F_c:=F_a+F_b$ and apply the the above process to $F_c$.
Then we get  a computation sequence $\{ F_c^{(0)}, \dots, 
 F_c^{(k_c)} \}$ (we denote $F_c^{(k)}:= 
{\pi_0}^*{\pi_0}_*(F_a+F_b)+\sum_{i=1}^n c_i^{(k)}{\pi_i}^*E_i$ 
for $0 \le k \le k_c$) for $F_a+F_b-\lceil K \rceil$. 
We define $\mathbf{a}^{(0)}$ and $\mathbf{b}^{(0)}$ as in $(2.2)$.

By definition and Claim $3$, $\mathbf{e}_j\mathbf{P}\mathbf{c}^{(0)}<0$ if and only if $\mathbf{e}_j\mathbf{P}(\mathbf{a}+\mathbf{b})=0$, $j \in \Lambda$, $a_j+b_j > 0$ and $a_i+b_i=0$ for every $i \in \Lambda$ such that $E_i \vartriangleright E_j$.
Therefore the condition $\mathbf{e}_j\mathbf{P}\mathbf{c}^{(0)}<0$ implies $\mathbf{e}_j\mathbf{P}\mathbf{a}^{(0)}<0$ and $\mathbf{e}_j\mathbf{P}\mathbf{b}^{(0)}<0$, unless $a_j$ or $b_j$ is zero.
If $c_j^{(1)}=c_j^{(0)}+1$, then we define $F_a^{(1)}:={\pi_0}^*{\pi_0}_*{F_a}+\sum_{i=1}^n a_i^{(1)}{\pi_i}^*E_i$ (resp. $F_b^{(1)}:={\pi_0}^*{\pi_0}_*F_b+\sum_{i=1}^n b_i^{(1)}{\pi_i}^*E_i$) as follows (cf. (2.3)):
If $a_j=0$ (resp. $b_j=0$), then set $F_a^{(1)}=F_a^{(0)}$ (resp. $F_b^{(1)}=F_b^{(0)}$).
If $a_j > 0$ (resp. $b_j >0$), then 
$$a_i^{(1)}=\left\{
\begin{array}{ll}
a_i^{(0)}+1 & \text{if $i=j$} \\
a_i^{(0)}-1 & \text{if $E_i \vartriangleright E_j$ and $a_i^{(0)}>0$} \\
a_i^{(0)} & \text{otherwise}
\end{array}
\right .$$
$$\left(
\text{resp. }
b_i^{(1)} =\left\{
\begin{array}{ll}
b_i^{(0)}+1& \text{if $i=j$} \\
b_i^{(0)}-1& \text{if $E_i \vartriangleright E_j$ and $b_i^{(0)}>0$} \\
b_i^{(0)} & \text{otherwise}
\end{array}
\right .
\right ).$$
We can inductively define a cycle $F_a^{(k)}$ (resp. $F_b^{(k)}$) for every $k=1, \dots, k_c$ as above.
When $a_j$ and $b_j$ are nonzero integers, $\mathbf{e}_j\mathbf{P}\mathbf{c}^{(k)}<0$ if and only if $\mathbf{e}_j\mathbf{P}\mathbf{a}^{(k)}<0$ and $\mathbf{e}_j\mathbf{P}\mathbf{b}^{(k)}<0$ for each $k=0,1,\dots, k_c$.
Hence $F_a^{(k_c)}$ (resp. $F_b^{(k_c)}$) is not necessarily $f$-anti-nef, but $\{F_a^{(1)}, \dots, F_a^{(k_c)}\}$ (resp. $\{F_b^{(1)}, \dots, F_b^{(k_c)} \})$ can be extended to a computation sequence $\{F_a^{(1)}, \dots, F_a^{(k_c)}, \dots, F_a^{(k_a)} \}$ (resp. $\{F_b^{(1)}, \dots, F_b^{(k_c)}, \dots, F_b^{(k_b)} \}$) for $F_a-\lceil K \rceil$ (resp. $F_b-\lceil K \rceil$).
Moreover, by Claim $3$, the triple $(a_i^{(k_c)}, b_i^{(k_c)}, c_i^{(k_c)})$ 
coincides with one of the following: 
$$(a_i, b_i, a_i+b_i), (a_i-1, b_i-1, a_i+b_i-1), (a_i-1, 0, a_i-1), (0, b_i-1, b_i-1).$$
Then we will show that $a_i^{(k_a)}+b_i^{(k_b)} \le c_i^{(k_c)}$.

If $(a_i^{(k_c)}, b_i^{(k_c)}, c_i^{(k_c)})=(a_i, b_i, a_i+b_i)$, then, by Claim $3$, we have $a_i^{(k_a)}+b_i^{(k_b)} \le a_i+b_i=c_i^{(k_c)}$.

If $(a_i^{(k_c)}, b_i^{(k_c)}, c_i^{(k_c)})=(a_i-1, 0, a_i-1)$, then $b_i=0$. Therefore we have $b_j=0$ for all $j$ such that $E_j > E_i$, because $F_b$ is $f$-anti-nef. Hence $b_i^{(k)}=b_i=0$ for all $k_c \le k \le k_b$ and $a_i^{(k)}=c_i^{(k_c)}=a_i-1$ for every $k_c \le k \le k_a$, in particular $a_i^{(k_a)}+b_i^{(k_b)}=a_i-1=c_i^{(k_c)}$. The case where $(a_i^{(k_c)}, b_i^{(k_c)}, c_i^{(k_c)})=(0, b_i-1, b_i-1)$ is similar.

Thus we suppose that $(a_i^{(k_c)}, b_i^{(k_c)}, c_i^{(k_c)})=(a_i-1, b_i-1, a_i+b_i-1)$. 
Then it suffices to prove that there exists no such $p,q \ge k_c$ as $a_i^{(p)}=a_i$ and $b_i^{(q)}=b_i$. Assume to the contrary that $a_i^{(p)}=a_i$ and $b_i^{(q)}=b_i$ for some $p,q \ge k_c$. 
Then among the index $j$'s for which $a_j^{(p)}=a_j$ and $b_j^{(q)}=b_j$ for some $p,q \ge k_c$, take $j$ such that $E_j$ is maximal with respcet to the relation $>$.

It does not occur that $a_j^{(k_c+1)}=a_j^{(k_c)}+1$ and $b_j^{(k_c+1)}=b_j^{(k_c)}+1$, because if $\mathbf{e}_j\mathbf{P}\mathbf{a}^{(k_c)}<0$ and $\mathbf{e}_j\mathbf{P}\mathbf{b}^{(k_c)}<0$, then $\mathbf{e}_j\mathbf{P}\mathbf{c}^{(k_c)}<0$, that is, $F_c^{(k_c)}$ is not $f$-anti-nef. 
Therefore if $\mathbf{e}_j\mathbf{P}\mathbf{a}^{(k_c)} < 0$, then $\mathbf{e}_j\mathbf{P}\mathbf{b}^{(k_c)} \ge 0$, which implies by Claim $3$ that $\mathbf{e}_j\mathbf{P}\mathbf{b}>0$ or there exists $E_h \vartriangleright E_j$ such that $(a_h^{(k_c)}, b_h^{(k_c)}, c_h^{(k_c)})=(0, b_h-1, b_h-1)$. Notice that, in the latter case, $b_h^{(k)}=b_h-1$ for every $k_b \ge k \ge k_c$. 

Thus, in both cases, $\mathbf{e}_j\mathbf{P}\mathbf{b}^{(k)} \ge 0$ for all $k_b \ge k \ge k_c$, namely $b_j^{(k)}=b_j-1$ for every $k_b \ge k \ge k_c$. 
Hence if $a_j^{(p)}=a_j$ and $b_j^{(q)}=b_j$ for some $p,q \ge k_c$, then $\mathbf{e}_j\mathbf{P}\mathbf{a}^{(k_c)} \ge 0$, $\mathbf{e}_j\mathbf{P}\mathbf{b}^{(k_c)} \ge 0$ and $\mathbf{e}_j\mathbf{P}\mathbf{a}=\mathbf{e}_j\mathbf{P}\mathbf{b}=0$. 
If there exists some $E_h \vartriangleright E_j$ such that $(a_h^{(k_c)}, b_h^{(k_c)}, c_h^{(k_c)})=(a_h-1, 0, a_h-1)$ or $(0, b_h-1, b_h-1)$, then it does not occur that $a_j^{(p)}=a_j$ and $b_j^{(q)}=b_j$ for some $p,q \ge k_c$.

Therefore there exists $E_h \vartriangleright E_j$ such that $(a_h^{(k_c)}, b_h^{(k_c)}, c_h^{(k_c)})=(a_h-1, b_h-1, a_h+b_h-1)$, and in order that $a_j^{(p)}=a_j$ and $b_j^{(q)}=b_j$ for $p,q \ge k_c$, it is necessary that $a_h^{(s)}=a_h$ and $b_h^{(t)}=b_h$ for some $p>s>k_c$ and $q>t>k_c$. 
However this contradicts the maximality of the index $j$.
Thus we have $a_i^{(k_a)}+b_i^{(k_b)} \le a_i+b_i-1=c_i^{(k_c)}$.
\end{proof}

\begin{rem}
$(1)$ If $A$ is a two-dimensional Gorenstein log terminal singularity, then the relative canonical divisor $K$ is an integral divisor. In this case, the anti-nef closure of $Z-K$ for an $f$-anti-nef cycle $Z$ can simply be described as follows:

$$\mathrm{an}_f(Z-K) = Z-K +\sum_{\scriptstyle {\pi_i}_*Z \cdot E_i=0 }{\pi_i}^{*}E_i .$$
By this formula, one can see that there are many cases in which the equality in the subadditivity theorem fails (see Example \ref{2-dim ex} $(1)$). 
On the other hand, in non-Gorenstein case, $\mathrm{an}_f(Z-\lceil K \rceil)$ is more complicated and an analog of the above formula 
$$\mathrm{an}_f(Z-\lceil K \rceil) = Z-\lceil K \rceil +\hspace{-1cm} \sum_{\scriptstyle  \lceil K_i \rceil = {f_i}^{*}\lceil K_{i-1} \rceil+E_i \atop \scriptstyle {\pi_i}_*Z \cdot E_i=0 }\hspace{-1cm}{\pi_i}^{*}E_i,$$
where $K_i$ is the relative canonical divisor of $f_i$, does not hold (see Example \ref{2-dim ex} $(2)$).

$(2)$ Multiplier ideals can also be defined for a divisor (refer to \cite{La} for details), and Demailly, Ein and Lazarsfeld \cite{DEL} also proved the subadditivity theorem for divisors under the assumption that $A$ is regular. Some readers may expect that the subadditivity theorem for divisors also holds if $A$ is a two-dimensional log terminal singularity, but this is not true. For example, let $A=\mC[[X,Y,Z]]/(XY-Z^2)$, $\mathcal{O}_{\Spec A}(-D_1)=(x,z)$ and $\mathcal{O}_{\Spec A}(-D_2)=(y,z)$, where $x,y,z$ are the images of $X,Y,Z$ in $A$, respectively. Then $\mathcal{J}(D_1)=(x,z)$ and $\mathcal{J}(D_2)=(y,z)$, but $\mathcal{J}(D_1+D_2)=(z)$. Hence $\mathcal{J}(D_1+D_2) \not\subseteq \mathcal{J}(D_1)\mathcal{J}(D_2)$.
\end{rem}

\begin{expl}\label{2-dim ex}
$(1)$ $($$A_1$-singularity$)$
Let $A=\mC[[X,Y,Z]]/(XY-Z^2)$ and $\mathfrak{a}=(x^2, y, z)$, where $x$, $y$, $z$ are the images of $X$, $Y$, $Z$ in $A$, respectively. Then
\begin{align*}
\mathcal{J}(\mathfrak{a})&=(x, y, z), \\
\mathcal{J}(\mathfrak{a}^2)&=(x^3, y^2, z^2, xy, yz, zx).
\end{align*}
Therefore $\mathcal{J}(\mathfrak{a}^2) \subsetneq \mathcal{J}(\mathfrak{a})^2$. Indeed, let $f_0:X_0 \to \Spec A$ be the minimal resolution with the irreducible exceptional curve $F$ and $f_1:X \to X_0$ the blowing-up of $X_0$ at a point on $F$. 
Let $f:X \to \Spec A$ be the composite morphism of $f_0$ and $f_1$.
Then the relative canonical divisor $K$ of $f$ is $E_1$ and the ideal $\mathfrak{a}$ corresponds to the $f$-anti-nef cycle $F_a=2E_1+F'$, where $E_1$ is the unique $f_1$-exceptional curve and $F'$ is the strict transform of $F$. Hence we have $\mathrm{an}_f(F_a-K)=E_1+F'$ and $\mathrm{an}_f(2F_a-K)=3E_1+2F'$. Thus $2 \ \mathrm{an}_f(F_a-K)<\mathrm{an}_f(2F_a-K)$.

$(2)$ Let $A=\mC[X^5, XY^3, X^2Y, Y^5] \subset \mC[X,Y]$. We denote by $F_1$ and $F_2$ the exceptional curves of the minimal resolution $f_0:X_0 \to \Spec A$ such that ${F_1}^2=-2$ and ${F_2}^2=-3$. 
Let $f_1:X_1 \to X_0$ be the blowing-up at the intersection of $F_1$ and $F_2$ with the $f_1$-exceptional curve $E_1$, and $f_2:X \to X_1$ be the blowing-up at the intersection of $E_1$ and the strict transform ${F_2}'$ of $F_2$ with the $f_2$-exceptional curve $E_2$. 
Let $f:X \to \Spec A$ be the composite morphism of $f_0, f_1$ and $f_2$. Then the graph of $f$-exceptional curves is the following.
$$
\begin{footnotesize}
\entrymodifiers={++[o][F-]}
\xymatrix@R=5pt{
 -3 \ar@{-}[r]
& -2 \ar@{-}[r]
& -1 \ar@{-}[r]
& -5 \\
*{{F_1}'} & *{{E_1}'} & *{E_2} &*{{F_2}'}}
\end{footnotesize}
$$
Fix an $f$-anti-nef cycle $Z={F_1}'+3{E_1}'+5E_2+{F_2}'$ on $X$. Since the relative canonical divisor $K$ of $f$ is $-\frac{1}{5}{F_1}'+\frac{2}{5}{E_1}'+E_2-\frac{2}{5}{F_2}'$, we have
\begin{align*}
\mathrm{an}_f(Z-\lceil K \rceil)&=\mathrm{an}_f({F_1}'+2{E_1}'+4E_2+{F_2}')={F_1}'+3{E_1}'+4E_2+{F_2} \\
&=Z-\lceil K \rceil+{E_1}'\\
&\ne Z-\lceil K \rceil +\hspace{-1cm} \sum_{\scriptstyle  \lceil K_i \rceil = {f_i}^{*}\lceil K_{i-1} \rceil+E_i \atop \scriptstyle {\pi_i}_*Z \cdot E_i=0 }\hspace{-1cm}{\pi_i}^{*}E_i.
\end{align*}
\end{expl}

On the other hand, in order that the strong subadditivity theorem 
holds, regularity is necessary.
\begin{prop}
Let $(A, \mathfrak{m})$ be a two-dimensional $\mQ$-Gorenstein normal local ring such that the residue field $A/\mathfrak{m}$ is algebraically closed.
If the strong subadditivity theorem holds, that is, 
$$\mathcal{J}(\mathfrak{a}^c \mathfrak{b}^d) \subseteq \mathcal{J}(\mathfrak{a}^c)  \mathcal{J}(\mathfrak{b}^d)$$ 
for any ideals $\mathfrak{a}, \mathfrak{b} \subseteq A$ and  any rational numbers $c, d>0$, then $A$ is regular.
In particular when $A$ is essentially of finite type over a field of characteristic zero, $A$ is regular if and only if the strong subadditivity theorem holds.
\end{prop}
\begin{proof}
Assume that $A$ is not regular. Let $f:X \to \Spec A$ be the minimal resolution, and then the exceptional locus $\mathrm{Exc}(f)$ of $f$ is not trivial. In order that the strong subadditivity theorem hold, by Theorem \ref{m-primary}, it is necessary that $A$ is a log terminal singularity.
\begin{enumerate}
\item the case where $\mathrm{Exc}(f)$ is irreducible.

\noindent Let $E$ be the unique irreducible $f$-exceptional
 curve.  Then $E^2=-k$ for some integer $k \geq 2$.
Let $g:Y \to X$ be the blowing-up at a point on the curve $E$ and $h:Y \to \Spec A$ the composite morphism of $f$ and $g$. 
We denote by $E_1$ the exceptional divisor of $g$ and by $E_2$ the strict transform of $E$. Then ${E_1}^2=-1$, ${E_2}^2=-k-1$ and the relative canonical divisor $K$ of $h$ is equal to $\frac{2}{k}E_1-\frac{k-2}{k}E_2$. Fix an $h$-anti-nef cycle $Z:=2(k+1)E_1+2E_2$ on $Y$, which is an $h$-anti-nef cycle on $Y$. Here, for every cycle $F$ on $Y$, we denote by $\mathrm{an}_h(F)$ the $h$-anti-nef closure of $F$ as in the proof of Theorem \ref{m-primary}. Then
\begin{align*}
& \mathrm{an}_{h}(\lfloor \frac{1}{k+1}Z-K \rfloor)=\mathrm{an}_{h}(E_1)=E_1+E_2, \\
& \mathrm{an}_{h}(\lfloor \frac{2}{k+1}Z-K \rfloor)=\mathrm{an}_{h}(3E_1+E_2)=3E_1+E_2.
\end{align*}
Hence $\mathcal{J}(I^{\frac{2}{k+1}}) \not\subseteq \mathcal{J}(I^{\frac{1}{k+1}})^2$, where $I=H^0(Y, \mathcal{O}_Y(-Z)) \subset R$. This implies that the strong subadditivity theorem does not hold on $A$.

\item the case where $\mathrm{Exc}(f)$ is reducible.

\noindent Let $Z_f$ be the fundamental cycle of $f$. Since $Z_f$ is reducible, we can take an $f$-anti-nef cycle $Z$ such that $Z_f \leq Z < nZ_f$ and $\lfloor \frac{1}{n}Z \rfloor \ne 0$ for some integer $n \geq 2$. We denote by $K_0$ the relative canonical divisor of $f$. Since $A$ is log terminal and $-K_0$ is an effective divisor,
\begin{align*}
& \mathrm{an}_f(\lfloor \frac{1}{n}Z - K_0 \rfloor) \ge \mathrm{an}_f(\lfloor \frac{1}{n}Z \rfloor) =Z_f, \\
& \mathrm{an}_f(\lfloor Z - K_0 \rfloor)=\mathrm{an}_f(Z)=Z.
\end{align*}
Therefore, denoting the ideal $I=H^0(Y, \mathcal{O}_Y(-Z)) \subset R$, we have $\mathcal{J}(I) \not\subseteq \mathcal{J}(I^{\frac{1}{n}})^n$. Thus the strong subadditivity theorem does not hold on $A$.
\end{enumerate}
\end{proof}

\begin{rem}
We believe that the strong subadditivity theorem also holds for any two-dimensional regular local ring which is not necessarily essentially of finite type over a field of characteristic zero, but we cannot prove this by an argument about anti-nef cycles such as the proof of Theorem \ref{m-primary}.
\end{rem}
\section{higher dimensional case}
In higher dimensional case, we have several counterexamples to Theorem \ref{m-primary}.

\begin{expl}\label{3-dim}
Let $A=\mC [[X,Y,Z,W]]/(X^2+Y^4+Z^4+W^5)$ and $\mathfrak{m}=(x,y,z,w)$, where $x$, $y$, $z$, $w$ are the images of $X$, $Y$, $Z$, $W$ in $A$, respectively. Then $A$ is a Gorenstein log terminal singularity, but not a terminal singularity. Therefore $\mathcal{J}(\mathfrak{m})=\mathfrak{m}$. If $\mathcal{J}(\mathfrak{m}^{k+l}) \subseteq \mathcal{J}(\mathfrak{m}^k) \mathcal{J}(\mathfrak{m}^l)$ holds for all integers $k,l> 0$ , then $\mathcal{J}(\mathfrak{m}^n)=\mathfrak{m}^n$ for any integer $n>0$, 
in particular $\mathfrak{m}^2$ should be integrally closed. However, $x \in \overline{\mathfrak{m}^2} \setminus \mathfrak{m}^2$, because $x^2 \in \mathfrak{m}^4$. This is a contradiction, hence the subadditivity theorem fails on $A$.
\end{expl}

Now we investigate the subadditivity theorem for monomial ideals. We expected that the subadditivity theorem for monomial ideals might hold on every toric ring, but unfortunately we found a counterexample on a three-dimensional toric ring.
\begin{expl}\label{toric}
Let $M=\{(x,y,z) \in \mZ^3 \mid 35x+28y+20z \equiv 0 \mod 41 \}$ be a lattice, $\sigma^{\vee}=(\mZ_{\ge 0})^{\oplus 3} \subset M \otimes_{\mZ} \mR$ a cone and $A=k[M\cap \sigma^{\vee}] \subset k[x,y,z]$ the cyclic quotient singularity of type $1/41(35, 28, 20)$. We consider the monomial ideal 
$$I=(x^{410}, y^{410}, z^{410}, x^{8}yz, x^4y^6z, x^4yz^8)\subset A.$$
Then we will prove that $\mathcal{J}(I^2) \nsubseteq \mathcal{J}(I)^2$. 

First we will show that $x^{10}y^{3}z^{7} \in \mathcal{J}(I^2)$. Let $P(I)$ be the Newton polygon of $I$, that is, the convex hull of 
$$\{(410,0,0), (0,410,0), (0,0,410), (8,1,1), (4,6,1), (4,1,8) \}$$
in $M \otimes_{\mZ} \mR=\mR^3$. Then note that for every positive integer $n$, by [HY, Theorem 4.8] and \cite{How}, $x^{a}y^bz^c \in \mathcal{J}(I^n)$ if and only if the point $(a+1, b+1, c+1)$ is contained in the interior $\mathrm{Int}(nP(I))$ of $nP(I)$. 
Since $(10+1, 3+1, 7+1)=\frac{245}{328}(8,1,1)+\frac{131}{328}(4,6,1)+\frac{281}{328}(4,1,8)$ and $\frac{245}{328}+\frac{131}{328}+\frac{281}{328}>2$, by the above characterization of multiplier ideals associated to a monomial ideal, $x^{10}y^{3}z^{7}$ is contained in $\mathcal{J}(I^2)$.

Next we will show that $x^{10}y^{3}z^{7}$ is not contained in $\mathcal{J}(I)^2$. If $x^{10}y^{3}z^{7} \in \mathcal{J}(I)^2$, then there exist lattice points $(p,q,r), (s,t,u) \in M$ such that $(p,q,r)+(s,t,u)=(10,3,7)$ and $(p+1,q+1,r+1), (s+1,t+1,u+1) \in \mathrm{Int}(P(I))$. 
Since the three points $(8,1,1), (4,6,1), (4,1,8)$ lie on the plane $35x+28y+20z=328$, by the condition that $(p+1,q+1,r+1), (s+1,t+1,u+1) \in \mathrm{Int}(P(I))$, the lattice points $(p,q,r)$ and $(s,t,u)$ must satisfy that $35p+28q+20r>328-(35+28+20)=245$ and $35s+28t+20u>245$. 
Moreover by the assumption $(p,q,r)+(s,t,u)=(10,3,7)$, we know that $(p,q,r)$ and $(s,t,u)$ are obliged to be $(8,1,1)$ and $(2,2,6)$. However $(2+1,2+1,6+1)$ is not contained in $\mathrm{Int}(P(I))$, because $(2+1,2+1,6+1)=-\frac{83}{328}(8,1,1)+\frac{131}{328}(4,6,1)+\frac{281}{328}(4,1,8)$.
This is a contradiction. Thus $x^{10}y^{3}z^{7} \notin \mathcal{J}(I)^2$.
\end{expl}
\begin{ques}
Let $A$ be a Gorenstein toric ring and $\mathfrak{a}, \mathfrak{b}$ be monomial ideals of $A$. Then
$$\mathcal{J}(\mathfrak{ab}) \subseteq \mathcal{J}(\mathfrak{a}) \mathcal{J}(\mathfrak{b})?$$
\end{ques}

\end{document}